\documentclass[graybox]{svmult}

\usepackage{graphicx}        
\usepackage[bottom]{footmisc}

\usepackage{newtxtext}       %
\usepackage{newtxmath}       

\usepackage{accents}
\usepackage{amscd}
\usepackage{cleveref}
\usepackage{mathtools}
\usepackage{stackrel}
\usepackage{url}

\DeclareMathOperator{\Expectation}{\mathbb E}
\DeclareMathOperator{\Grad}{grad}
\DeclareMathOperator{\Maxexp}{\mathcal E}

\newcommand{\KL}[2]{\operatorname{D}\left(#1\,\Vert#2\right)}
\newcommand{\anytransport}[2]{\prescript{\cdot}{} {\mathbb U} _ {#1} ^ {#2}}
\newcommand{\condexpat}[3]{\Expectation_{#1}\left(#2 \middle| #3\right)}
\newcommand{\derivby}[1]{\frac{d}{d#1}}
\newcommand{\etransport}[2]{\prescript{\text{e}}{} {\mathbb U} _ {#1} ^ {#2}}
\newcommand{\euler}{\mathrm{e}}
\newcommand{\expbundleat}[1]{S\maxexpat{#1}}
\newcommand{\expectat}[2]{\Expectation_{#1}\left[#2\right]}
\newcommand{\expfiberat}[2]{S_{#1}\maxexpat{#2}}
\newcommand{\expof}[1]{\exp\left(#1\right)}
\newcommand{\maxexpat}[1]{\Maxexp\left(#1\right)}
\newcommand{\mtransport}[2]{\prescript{\text{m}}{} {\mathbb U} _ {#1} ^ {#2}}
\newcommand{\predual}[1]{\prescript{*}{}{#1}}
\newcommand{\reals}{\mathbb{R}}
\newcommand{\scalarat}[3]{\left\langle#2,#3\right\rangle_{#1}}
\newcommand{\setof}[2]{\left\{#1 \, \middle| \, #2 \right\}}
\newcommand{\velocity}[1]{\accentset{\star}{#1}}

\begin{document}

\title*{Information geometry of Bayes computations}
\titlerunning{IG of Bayes}
\author{Giovanni Pistone}

\institute{Giovanni Pistone \at Nuovo SEFIR, C/O COWORLD Segrate IT \email{giovanni.pistone@carloalberto.org}}

\maketitle

\abstract{Amari's Information Geometry is a dually affine formalism for parametric probability models. The literature proposes various nonparametric functional versions. Our approach uses classical Weyl's axioms so that the affine velocity of a one-parameter statistical model equals the classical Fisher's score. In the present note, we first offer a concise review of the notion of a statistical bundle as a set of couples of probability densities and Fisher's scores. Then, we show how the nonparametric dually affine setup deals with the basic Bayes and Kullback-Leibler divergence computations.}

\section{Dually Affine Information Geometry: the Statistical Bundle Formalism}

We review some introductory material from our tutorial \cite{chirco|pistone:2022}. See there a complete list of references. If $(\Omega,\mathcal F, \mu)$ is a probability space and $B$ a functional Banach space of random variables, a \emph{maximal exponential model} $\maxexpat{\mu}$ is a set of probability densities of the form $\expof{u - K_p(u)} \cdot p$, $p \in \maxexpat \mu$, $u \in \mathcal S \subset B$, $\expectat p u = 0$. We assume $\mathcal S$ is the largest open, convex domain, and it is such that an open exponential arc connects all couples of densities. Any probability density in the model can tell the role of $p$. There are many options for the Banach space $B$. The most straightforward option is to assume $\Omega$ as a finite set or a bounded real domain and take $B$ as the space of real vectors or continuous real functions. The relevant literature knows many variants for the space $B$, for example, Orlicz spaces \cite{pistone|sempi:95,pistone:2013gsi}, Hilbert spaces \cite{egozcue|diaz-barrero|pawlowsky-glahn:2006,talska|menafoglio|hron|egozcue|palarea-albaladejo:2019}, spaces of measures \cite{ay|jost|le|schwachhofer:2017igbook}, spaces of infinitely differentiable functions \cite{khesin|misoliek|modin:2024}.

On the maximal exponential model $\maxexpat \mu$ we define a special vector bundle, \emph{statistical bundle},
\begin{equation*}
    \expbundleat \mu = \setof{(q,v)}{q \in \maxexpat \mu, v \in B, \expectat q v = 0} \ .
\end{equation*}

The original motivation for this formalism is Fisher's score. See, for example, the textbook \cite[Ch.~4]{efron|hastie:2016}. If $\theta \mapsto q(\theta)$ is a 1-dimensional smooth curve in $\maxexpat \mu$, then
\begin{equation*}
    \theta \mapsto \left(q(\theta),\derivby \theta \log q(\theta)\right) \in \expbundleat \mu \ .
\end{equation*}
We interpret Fisher's score to be the \emph{velocity} of the curve so that the fibers $\expfiberat{q}{\mu}$ are the space of velocities, i.e., the tangent spaces. The statistical bundle has the same function as the state space in Mechanics. See, for example, the textbook \cite[Ch.~IV]{arnold:1989}. 

On each fiber $\expfiberat{q}{\mu}$ of the statistical bundle has the covariance bilinear form, which is the restriction of a duality pairing between two Banach spaces of random variables $\predual B$ and $B$, 
\begin{equation*}
 \predual{\expfiberat{q}{\mu}} \times \expfiberat{q}{B} \ni (v,w) \mapsto \scalarat q v w = \expectat q {vw} \ .   
\end{equation*}
In our theory, the basic structure is the affine space defined below, not the Hilbert space $L^2(p \cdot \mu)$. The spaces in duality are equal in the straightforward case we consider here, $\predual B = B$, but it is nevertheless helpful to keep them distinct.

We have two \emph{parallel transports}
\begin{gather}
    \etransport pq \colon \expfiberat p {\mu} \ni v \mapsto v - \expectat q v \in \expfiberat q \mu \\
    \mtransport pq \colon \predual{\expfiberat p {\mu}} \ni w \mapsto \frac p q w \in \predual{\expfiberat q {\mu}}   
\end{gather}
It is easy to verify that both parallel transports form a co-cycle and are dual to each other:
\begin{equation*}
 \anytransport q r \anytransport p q = \anytransport p r \ , \quad \anytransport pp = \operatorname{Id} \ , \quad   \scalarat q {\mtransport p q w} v = \scalarat p w {\etransport q p v} \ .
\end{equation*}

The dually affine structure follows from the definition of two \emph{affine displacements} associating to each couple of points a vector. If the first point is the frame's origin, the displacement becomes an affine chart. We define two atlas of charts: for all $p,q \in \maxexpat \mu$, 
\begin{gather}
    s_p(q) = \log \frac q p - \expectat p {\log \frac q p} \in \expfiberat p {\mu} \quad \text{is an exponential chart and}\\
    \eta_p(q) = \frac q p - 1 \in \predual{\expfiberat p {\mu}} \quad \text{is a mixture chart.} \label{eq:mixture-chart}  
\end{gather}
The names exponential and mixture come from the form of the inverse charts,
\begin{equation*}
s_p^{-1}(v) = \euler^{v - K_p(v)}\cdot p \ , \quad \eta_p^{-1}(w) = (1+w) \cdot p \ .
\end{equation*}

The \emph{Weyl's axioms} hold in the generalized form
\begin{gather*}
    s_p(q) + \etransport q p s_q(r) = s_p(r) \ , \\
    \eta_p(q) + \mtransport q p \eta_q(r) = \eta_p(r) \ ,
\end{gather*}
hence, all the machinery of the calculus of affine spaces holds for the statistical bundles. In particular, the velocity of a curve $t \mapsto q(t) \in \maxexpat {B}$ expressed in the \emph{moving frame} equals the Fisher's score,
\begin{equation*}
   \left. \derivby t s_p(q(t)) \right|_{p=q(t)} =  \left. \derivby t \eta_p(q(t)) \right|_{p=q(t)} = \derivby t \log q(t) = \velocity q(t) \ .
\end{equation*}

\section{The gradient of the KL Divergence} \label{sec:KL-gradient}

For all $p,q \in \maxexpat \mu$ we can write
\begin{equation} \label{eq:structural-equation}
    q = \euler^{s_p(q) - \KL p q} \cdot p \ , 
\end{equation}
where the additive normalizing constant is the Kullback-Leibler (KL) divergence
\begin{equation*}
  \maxexpat \mu \times \maxexpat \mu \ni (q,r) \mapsto \KL q r = \expectat q {\log \frac q r} \ .
\end{equation*}
Cf. \cite[Ch.~5]{amari:2016}.

The total natural gradient of the KL divergence is a couple of sections $(\Grad_1 \KL \cdot \cdot,\Grad_2 \KL \cdot \cdot)$ of the statistical bundles $\expbundleat \mu$ and $\predual{\expbundleat{\mu}}$, respectively, such that for all couple of smooth curves $t \mapsto (q(t),r(t)) \in \maxexpat \mu \times \maxexpat \mu$ it holds
\begin{multline*} 
    \derivby t \KL {q(t)}{r(t)} = \\ \scalarat {q(t)} {\velocity q(t)} {\Grad_1 \KL {q(t)}{r(t)}} +
    \scalarat {r(t)} {\Grad_2 \KL {q(t)}{r(t)}} {\velocity r(t)} \ .
    \end{multline*}
    See a review of the definition of the natural gradient we use here in \cite[\S~12.1.2]{amari:2016} and \cite{chirco|pistone:2022}.
    
A direct computation in the exponential and mixture charts, respectively, provides the neat result
\begin{align} 
    \Grad_1 \KL q r &= - s_{q}(r) \ , \label{eq:KL-grad-1}\\
    \Grad_2 \KL q r &= - \eta_{r}(q) \ . \label{eq:KL-grad-2}
\end{align}
A detailed derivation will appear in \cite{pistone:2025-MDPI-Stats}. Because of the special relation between the KL divergence and the affine charts, the KL divergence has a special status of \emph{natural divergence} in the context of the dually-affine IG.

As a first example, assume $q = M(\theta)$ and $r = N(\theta)$, with $\theta \in \Theta \subset \reals^d$ are parametric models with a common parameter,
\begin{multline*}
  \nabla \KL {M(\theta)} {N(\theta)} = \\ - \scalarat {M(\theta)} {\nabla \log M(\theta)} {s_{M(\theta)}(N(\theta))} - \scalarat {N(\theta)} {\eta_{N(\theta)}(M(\theta))} {\nabla \log N(\theta)} = \\
  \int \left(\log \frac {N(\theta)} {M(\theta)} \nabla M(\theta) + (N(\theta) - M(\theta)) \nabla \log N(\theta)\right) d\mu \ .
\end{multline*}

\section{Bayes computations}

We first review the derivation of a function $f$ between two maximal exponential models. This section uses the mixture charts \eqref{eq:mixture-chart}. The expressions of $f$ and its derivative $df$ in the charts centered, respectively, at $p_1$ and $p_2$, are, 
\begin{equation*}
  \begin{array}{lcr}
    \begin{CD}
        \maxexpat {\mu_1} @>f>> \maxexpat {\mu_2} \\
        @A\eta_{p_1}^{-1}AA @VV\eta_{p_2}V \\
        \expfiberat {p_1} {\mu_1} @>>f_{p_1,p_2} > \expfiberat {p_2} {\mu_2}
      \end{CD}
    & \ \text{and} \ & 
      \begin{CD}
      \expfiberat {q_1} {\mu_1} @>df(q_1)>> \expfiberat {f(q_1)} {\mu_2} \\
      @A\mtransport{p_1}{q_1}AA @VV\mtransport{f(q_1)}{p_2}V \\
      \expfiberat {p_1} {\mu_1} @>>df_{p_1,p_2}(\eta_{p_1}(q_1))> \expfiberat{p_2}{\mu_2} 
    \end{CD}
  \end{array}
\end{equation*}
The computation of the derivative from its expression is
\begin{equation}
  df(q)[\velocity q] = \mtransport {p_2}{q} df_{p_1,p_2}(\eta_{p_1}(q)))[\mtransport {q} {p_1} \velocity q] \label{eq:bundle-derivative} \ .
\end{equation}

Inspired by \cite[\S~11.5]{amari:2016}, we consider those mappings from probability densities on a product space $(\Omega_1 \times \Omega_2,\mu_1 \otimes \mu_2)$ to each of the margins which arise in Bayesian computations. $\Omega_1$ is the sample space, and $\Omega_2$ is the parameters' space. However, Bayesian arguments are just one of the many applications of the information geometry of the product space, with others being functional ANOVA and transport theory; see \cite{pistone:2021-gsi}.

\subsection{Marginalization}
\label{sec:marginalization}

Consider first the marginalization $f = X_\#$, where $X \colon \Omega_1 \times \Omega_2 \to \Omega_1$ is the projection on the first factor, 
\begin{equation*}
 f \colon \maxexpat{\mu_1 \otimes \mu_2} \ni q \mapsto q_1 = \int q(\cdot,z) \, \mu_2(dz) \in \maxexpat{\mu_1} \ . 
\end{equation*}
The expession of $f$ in the charts centered at $p_1 \otimes p_2$ and $p_1$ is linear,
\begin{multline*}
f_{p_1\otimes p_2,p_1} \colon v \overset{\eta_{p_1 \otimes p_2}^{-1}}\longmapsto (1+v) \cdot p_1 \otimes p_2 \overset{f}\longmapsto \\ \int (1+v(\cdot,z)) \, p_1(\cdot) \, p_2(z) \, \mu_2(dz) \overset{\eta_{p_1}}\longmapsto \int v(\cdot,z) \, p_2(z) \, \mu_2(dz) \ .
\end{multline*}
Now the derivative of the marginalization function $f$ follows from \cref{eq:bundle-derivative},
\begin{multline} \label{eq:marginalization-derivative}
  df(q)[\velocity q] = \mtransport {p_1} {q_1} df_{p_1\otimes p_2,p_1} (\eta_{p_1 \otimes p_2}(q)) [\mtransport q {p_1 \otimes p_2} \velocity q] = \\  \frac{p_1}{q_1} \int \frac{q(\cdot,z)}{p_1(\cdot) p_2(z)} \velocity q(\cdot,z) \, p_2(z) \, \mu_2(dz) =
\\ \left(x \mapsto \int \velocity q(x,z) q_{2|1}(z|x) \, \mu_2(dz)\right) = \condexpat q {\velocity q}{X} \ .
\end{multline}
The interpretation of the conditional expectation as the derivative of the marginalization provides valuable insight. 

\subsection{Conditioning}
\label{sec:conditioning}

We now turn to the conditioning function. For each $x \in \Omega_1$, the mapping
\begin{equation*}
f_x \colon  \maxexpat {\mu_1 \otimes \mu_2} \ni q_{12} \mapsto \left(y \mapsto q_{2|1}(y|x) = \frac{q_{12}(x,y)}{\int q_{12}(x,y) \mu_1(dx)} \right) \in \maxexpat {\mu_2}
\end{equation*}
is well-defined in our setup because the densities are defined everywhere. We compute the derivative in the mixture chart with \cref{eq:bundle-derivative} origin $p_1 \otimes p_2$ in the joint space and $p_2$ on the margin. The expression $F_x$ of $f_x$ is,
\begin{multline*}
 F_x \colon  v \overset{\eta^{-1}_{p_1 \otimes p_2}}{\longmapsto} (1+v) \cdot p_1 \otimes p_2 \overset{f_x}{\longmapsto}  \\
  \frac{(1+v(x,\cdot))p_2(\cdot)}{\int (1 + v(x,z)) \, p_2(z) \,\mu_2(dz)} \overset{\eta_{p_2}}{\longmapsto} \\
   \frac{v(x,\cdot) - \int v(x,z) \, p_2(z) \, \mu_2(dz)}{1 + \int v(x,z) \, p_2(z) \, \mu_2(dz)} \ .
 \end{multline*}
 
The derivative of $F_x$ at $v$ in the direction $h$ is
\begin{equation*}
dF_x(v)[h] = \frac{h(x,\cdot) - \int h(x,z) \, p_2(z) \, \mu_2(dz) - F_x(v) \int h(x,z) \, p_2(z) \, \mu_2(dz)}{1 + \int v(x,z) \, p_2(z) \, \mu_2(dz)} \ .
\end{equation*}
As $1+v = q_{12}/p_1 \otimes p_2$, we have $1 + \int v(x,z) \, p_2(z) \, \mu_2(dz) = q_1(x)/p_1(x)$, so that
$F_x(\eta_{P_1 \otimes p_2}(q_{12})) = q_{2|1}(\cdot|x)/p_2(\cdot) - 1 = \eta_{p_2}(q_{2|1}(\cdot|x))$, hence
\begin{equation*}
  dF_x(\eta_{p_1 \otimes p_2}(q_{12}))[h] = \frac{p_1(x)}{q_1(x)}\left(h(x,\cdot) - \frac{q_{2|1}(\cdot|x)}{p_2(\cdot)} \int h(x,z) \, p_2(z) \, \mu_2(dz) \right) \ .
\end{equation*}
Now we plug-in $h = \mtransport {q_{12}} {p_1 \otimes p_2} \velocity q_{12} = q_{12} \velocity q_{12} / p_1 \otimes p_2$, to get
\begin{equation*}
  dF_x(\eta_{p_1 \otimes p_2}(q_{12}))\left[\mtransport {q_{12}} {p_1 \otimes p_2} \velocity q_{12}\right] = \frac {q_{2|1}(\cdot|x)}{p_2(\cdot)} \left(\velocity q_{12}(x,\cdot) - \int \velocity q_{12}(x,z) \, q_{2|1}(z|x) \, \mu_2(dz)\right) \ .
\end{equation*}
In conclusion, the bundle derivative is
\begin{equation} \label{eq:conditional-derivative}
  df_x(q)[\velocity q_{12}] = \velocity q_{12}(x,\cdot) - \int \velocity q_{12}(x,z) \, q_{2|1}(z|x) \, \mu_2(dz) \ . 
\end{equation}

\subsection{Exponential decomposition}

Other results follow using the exponential atlas and the structural equation \cref{eq:structural-equation}. We have
\begin{align*}
    &q_{12}(x,y) = \euler^{u_{12}(x,y) - \KL {p_1 \otimes p_2} {q_12}} \cdot p_1 \otimes p_2 \ ,  &u_{12} = s_{p_1 \otimes p_2}(q_{12}) \ ; \\
    &q_1(x) = \euler^{u_1(x) - \KL {p_1} {q_1}} \cdot p_1 \ ,  &u_1 = s_{p_1}(q_1) \ ; \\
    &q_{2|1}(y|x) = \euler^{u_{2|1}(y|x) - \KL {p_2} {q_{2|1}(\cdot|x)}} \cdot p_2 \ ,  &u_{2|1}(\cdot|x) = s_{p_2}(q_{2|1}(\cdot|x) \ .
\end{align*}
As $q_{12}(x,y) = q_1(x) q_{2|1}(y|x)$, from the zero expectations follows
\begin{gather*}
    u_{12}(x,y) = u_1(x) + u_{2|1}(y|x) - \left(\KL {p_2} {q_{2|1}(\cdot|x)} - \int \KL {p_2} {q_{2|1}(\cdot|x)} \, \mu_1(dx)\right) \ , \\
    \KL {p_1 \otimes p_2} {q_{12}} = \KL {p_1} {q_1} + \int p_1(x) \KL {p_2} {q_{2|1}(\cdot|x)} \, \mu_1(dx) \ .
\end{gather*}

\subsection{Examples and conclusion}
\label{sec:examples}

As a first application, consider the marginalization of an exponential family \cite[Ch.~5]{efron|hastie:2016}. In our setup, an exponential family is a parametric family of the form
\begin{equation} \label{eq:exponential-family}
 \theta \mapsto G(\theta) = \expof{\theta \cdot T - \psi(\theta)} \cdot p_1 \otimes p_2 \in \maxexpat {\mu_1 \otimes \mu_2} \ ,
\end{equation}
with $\theta \in \Theta$, a convex open subset of $\reals^d$, and we assume, without restriction of generality, $T_j \in \expfiberat {p_1 \otimes p_2} {\mu_1 \otimes \mu_2}$, especially $\expectat {p_1 \otimes p_2} {T_j} = 0$. The normalizing constant equals the KL-divergence
\begin{equation*}
  \psi(\theta) = \log \int \euler^{\theta \cdot T} \, p_1 \otimes p_2 \, d (\mu_1 \otimes \mu_2) = \KL {p_1 \otimes p_2} {G(\theta)}
\end{equation*}
and $\nabla \psi(\theta) = \expectat {G(\theta)} T$.

The expression of the exponential model in the mixture chart is
\begin{equation*}
  G_{p_1 \otimes p_2} \colon\theta \mapsto \eta_{p_1 \otimes p_2}(G(\theta)) = \expof{\theta \cdot T - \psi(\theta)} - 1 \ ,
\end{equation*}
whose derivative in the direction $\dot \theta$ is
\begin{equation} \label{eq:exp-model-derivative}
  dG_{p_1 \otimes p_2}(\theta)[\dot \theta] = \frac{G(\theta)}{p_1 \otimes p_2} (T - \expectat {G(\theta)} T) \cdot \dot \theta \ .
\end{equation}
Hence the velocity of $\theta \mapsto G(\theta)$ and the velocity of the margin $\theta \mapsto G_1(\theta) = f \circ G(\theta)$ are, respectively,
\begin{gather}
  dG(\theta)[\dot \theta] = (T - \expectat {G(\theta)} T) \cdot \dot \theta \quad \text{(from \cref{eq:exp-model-derivative})}\ , \label{eq:model-derivative}\\ dG_1(\theta) = \condexpat {G(\theta)} {T - \expectat {G(\theta)} T}{X} \cdot \dot \theta \quad \text{(from \cref{eq:bundle-derivative})}\label{eq:model-derivative-1} \ .
\end{gather}

If $r_1$ is a generic margin, let us compute the gradient of the parameterized KL-divergences. From \cref{eq:model-derivative-1} with \cref{eq:KL-grad-2} and \cref{eq:KL-grad-1}, respectively, we find  
  \begin{gather*}
    \derivby t \KL{r_1}{G_1(\theta(t))} = 
- \left(\int \condexpat {G(\theta(t))} {T - \expectat {G(\theta(t))} T} {X} r_1 \, d\mu_1\right)\cdot \dot \theta(t)  \ , \\
\derivby t \KL{G_1(\theta(t))}{r_1} = 
- \left(\int \condexpat {G(\theta(t))} {T - \expectat {G(\theta(t))} T} {X} \log \frac {r_1}{G_1(\theta(t))} \ , d\mu_1\right) \cdot \dot \theta(t) \ .
\end{gather*}

Let us consider the conditioning of the exponential family of \cref{eq:exponential-family}. Given $x$, the parametric family is
\begin{equation*}
\theta \mapsto G_{2|1}(\cdot|x;\theta) = \frac{G(x,\cdot;\theta)}{G_1(x;\theta)} \ .
\end{equation*}
Along a parametric curve $t \mapsto \theta(t)$ the velocity of the exponential family \cref{eq:exponential-family} is
\begin{equation*}
  \velocity G(\theta(t)) = (T - \expectat {G(\theta(t))} T) \cdot \dot \theta(t) \ ,
\end{equation*}
while the velocity of the conditioned density follows from \cref{eq:conditional-derivative} 
\begin{equation*}
  \velocity G_{2|1}(\cdot|x;\theta(t)) = \left(T(x,\cdot) - \int T(x,z) G_{2|1}(z|x;\theta) \, \mu_2(dz)\right) \cdot \dot \theta(t) \ .
\end{equation*}

\noindent\textbf{Final remark} As a conclusion, we point out that the introduction of the statistical bundle provides a natural setting for the nonparametric dually affine Information Geometry of \cite{amari|nagaoka:2000} and \cite{amari:2016}. It is apparent in the computation of the total gradient of the Kullback-Leibler divergence, in the computation of the differential forms of the Bayes formula,  and in the example of the computation of gradients along an exponential model. We suggest that this formalism is helpful in the computation made in applied statistical topics such as Variational Bayes. Applying the differential equations we have derived depends on choosing an exponential family with suitable special characters. The literature on such applications is too extensive to present in a short conference note like this one and will be the object of a longer paper currently in progress.
  
\noindent\textbf{Acknowledgments} De Castro Stratistics, Collegio Carlo Alberto, Italy partially supported the production of this paper. The author is a member of the AI\&ML\&MAT UMI interest group.

\bibliographystyle{spmpsci}

\begin{thebibliography}{10}
\providecommand{\url}[1]{{#1}}
\providecommand{\urlprefix}{URL }
\expandafter\ifx\csname urlstyle\endcsname\relax
  \providecommand{\doi}[1]{DOI~\discretionary{}{}{}#1}\else
  \providecommand{\doi}{DOI~\discretionary{}{}{}\begingroup
  \urlstyle{rm}\Url}\fi

\bibitem{amari|nagaoka:2000}
Amari, S., Nagaoka, H.: Methods of information geometry.
\newblock American Mathematical Society (2000).
\newblock Translated from the 1993 Japanese original by Daishi Harada

\bibitem{amari:2016}
Amari, S.i.: Information geometry and its applications, \emph{Applied
  Mathematical Sciences}, vol. 194.
\newblock Springer, [Tokyo] (2016).
\newblock \urlprefix\url{https://doi.org/10.1007/978-4-431-55978-8}

\bibitem{arnold:1989}
Arnold, V.I.: Mathematical methods of classical mechanics, \emph{Graduate Texts
  in Mathematics}, vol.~60.
\newblock Springer-Verlag, New York (1989).
\newblock Translated from the 1974 Russian original by K. Vogtmann and A.
  Weinstein, Corrected reprint of the second (1989) edition

\bibitem{ay|jost|le|schwachhofer:2017igbook}
Ay, N., Jost, J., L\^{e}, H.V., Schwachh\"{o}fer, L.: Information geometry,
  \emph{Ergebnisse der Mathematik und ihrer Grenzgebiete. 3. Folge. A Series of
  Modern Surveys in Mathematics [Results in Mathematics and Related Areas. 3rd
  Series. A Series of Modern Surveys in Mathematics]}, vol.~64.
\newblock Springer, Cham (2017).
\newblock \doi{10.1007/978-3-319-56478-4}.
\newblock \urlprefix\url{https://doi.org/10.1007/978-3-319-56478-4}

\bibitem{chirco|pistone:2022}
Chirco, G., Pistone, G.: Dually affine {I}nformation {G}eometry modeled on a
  {B}anach space (2022).
\newblock \doi{10.48550/ARXIV.2204.00917}.
\newblock \urlprefix\url{https://arxiv.org/abs/2204.00917}.
\newblock ArXiv:2204.00917

\bibitem{efron|hastie:2016}
Efron, B., Hastie, T.: Computer age statistical inference, \emph{Institute of
  Mathematical Statistics (IMS) Monographs}, vol.~5.
\newblock Cambridge University Press, New York (2016).
\newblock \urlprefix\url{https://doi.org/10.1017/CBO9781316576533}.
\newblock Algorithms, evidence, and data science

\bibitem{egozcue|diaz-barrero|pawlowsky-glahn:2006}
Egozcue, J.J., D\'iaz–Barrero, J.L., Pawlowsky–Glahn, V.: Hilbert space of
  probability density functions based on {A}itchison geometry.
\newblock Acta Mathematica Sinica, English Series \textbf{22}(4), 1175–1182
  (2006).
\newblock \doi{10.1007/s10114-005-0678-2}.
\newblock \urlprefix\url{http://dx.doi.org/10.1007/s10114-005-0678-2}

\bibitem{khesin|misoliek|modin:2024}
Khesin, B., Misio{\l}ek, G., Modin, K.: Information geometry of diffeomorphism
  groups (2024).
\newblock \urlprefix\url{https://arxiv.org/abs/2411.03265}

\bibitem{pistone:2013gsi}
Pistone, G.: Nonparametric information geometry.
\newblock In: F.~Nielsen, F.~Barbaresco (eds.) Geometric science of
  information, \emph{Lecture Notes in Comput. Sci.}, vol. 8085, pp. 5--36.
  Springer, Heidelberg (2013).
\newblock First International Conference, GSI 2013 Paris, France, August 28-30,
  2013 Proceedings

\bibitem{pistone:2021-gsi}
Pistone, G.: Statistical bundle of the transport model.
\newblock In: F.~Nielsen, F.~Barbaresco (eds.) Geometric Science of
  Information, pp. 752--759. Springer International Publishing, Cham (2021)

\bibitem{pistone:2025-MDPI-Stats}
Pistone, G.: Constrained minima of the {K}ullback-{L}eibler divergence (2025).
\newblock In progress

\bibitem{pistone|sempi:95}
Pistone, G., Sempi, C.: An infinite-dimensional geometric structure on the
  space of all the probability measures equivalent to a given one.
\newblock Ann. Statist. \textbf{23}(5), 1543--1561 (1995)

\bibitem{talska|menafoglio|hron|egozcue|palarea-albaladejo:2019}
Talska, R., Menafoglio, A., Hron, K., Egozcue, J.J., Palarea-Albaladejo, J.:
  Changing reference measure in bayes spaces with applications to functional
  data analysis (2019).
\newblock \urlprefix\url{https://arxiv.org/abs/1912.08003}

\end{thebibliography}

\end{document}